\documentclass[a4paper,12pt]{article}

\usepackage{graphicx}
\usepackage{latexsym}
\usepackage{amssymb}
\usepackage{amsmath}
\usepackage{layout}
\usepackage{enumerate}
\usepackage{latexsym}
\usepackage{amssymb}
\usepackage{amsmath}
\usepackage[latin1]{inputenc}
\usepackage{layout}

\newtheorem{teo}{Theorem}[section]
\newtheorem{lema}[teo]{Lemma}

\newtheorem{prop}[teo]{Proposition}
\newtheorem{obs}[teo]{Remark}

\newtheorem{hypothesis}[teo]{Hypothesis}

\newcommand{\R}{\mathbb{R}}
\newcommand{\N}{\mathbb{N}}

\newcommand{\al}{\alpha}
\newcommand{\Om}{\Omega}

\newcommand{\Del}{\Delta}
\newcommand{\ep}{\varepsilon}
\newcommand{\lam}{\lambda}

\newcommand{\gam}{\gamma}

\newcommand{\LL}{\mathcal{L}}

\newcommand{\C}{\mathcal{C}}

\newcommand{\tf}{\mathcal{F}}
\newcommand{\beq}{\begin{equation}}
\newcommand{\eeq}{\end{equation}}
\newcommand{\intdob}{\int_0^T\int_0^1}

\newenvironment{dem}{\text\bf {\it Proof:}}{}

\def\bistackrel#1#2#3{\mathrel{\mathop{#3}\limits^{#1}_{#2}}}
\begin{document}

\title{Weak convergence for the stochastic heat equation driven by Gaussian white noise}
\author{ Xavier Bardina, Maria Jolis and Llu\'{\i}s Quer-Sardanyons 
\thanks{The three authors are supported by the grant MEC-FEDER Ref. MTM2006-06427 from the   Direcci\'on General de Investigaci\'on, Ministerio de Educaci\'on y Ciencia, Spain.}\\
\small Departament de Matem\`atiques \\ 
\small Universitat Aut\`onoma de Barcelona \\
\small 08193 Bellaterra (Barcelona) Spain\\ 
\small bardina@mat.uab.cat; mjolis@mat.uab.cat; quer@mat.uab.cat
}

\maketitle

\begin{abstract}
In this paper, we consider a quasi-linear stochastic heat equation on $[0,1]$,
with Dirichlet boundary conditions and controlled by the space-time
white noise. We formally replace the random perturbation by 
a family of noisy inputs depending on a parameter
$n\in \N$ such that approximate the white noise in some sense. Then, we
provide sufficient conditions ensuring that the
real-valued {\it mild} solution of the SPDE perturbed by this family of noises converges 
in law, in the space $\C([0,T]\times [0,1])$ of
continuous functions, to the solution of the white noise driven SPDE. 
Making use of a suitable continuous functional of the stochastic
convolution term, we show that it suffices to tackle the linear
problem. For this, we prove that the corresponding family of laws is
tight and we identify the limit law by showing the convergence of
the finite dimensional distributions. We have also considered two
particular families of noises to that our result applies. The first
one involves a Poisson process in the plane and has been motivated
by a one-dimensional result of Stroock, which states that the family
of processes $n \int_0^t (-1)^{N(n^2 s)} ds$, where $N$ is a
standard Poisson process, converges in law to a Brownian motion. The
second one is constructed in terms of the kernels associated to the
extension of Donsker's theorem to the plane.

\end{abstract}

\bigskip



\pagebreak

\section{Introduction}

In the almost last three decades, there have been enormous advances
in the study of random field solutions to stochastic partial
differential equations (SPDEs) driven by general Brownian noises.
The starting point of this theory was the seminal work by Walsh
\cite{walsh}, and most of the research developed thereafter has been
mainly focused on the analysis of heat and wave equations perturbed
by Gaussian white noises in time with a fairly general spatial
correlation (see, for instance, \cite{bgp,cn,dc,dalang,ms}). Notice
also that some effort has been made to deal with SPDEs driven by
fractional type noises (see, for instance, \cite{glt,mn,qs,ttv}).

Indeed, the motivation to consider  these type of models in the
above mentioned references has sometimes put together theoretical
mathematical aspects and applications to some real situations. Let
us mention that, for instance, different type of SPDEs provide
suitable models in the study of growth population, some climate and
oceanographical phenomenons, or some applications to mathematical
finance (see \cite{dawson}, \cite{imkeller}, \cite{boris},
\cite{bjork}, respectively).

However, real noisy inputs are only approximately white and
Gaussian, and what one usually does is to justify somehow that one
can approximate the randomness acting on the system by a Gaussian
white noise. This fact has been illustrated by Walsh in
\cite{walsh-neural}, where a parabolic SPDE has been considered in
order to model a discontinuous neurophysiological phenomenon. The
noise considered in this article is determined by a Poisson point
process and the author shows that, whenever the number of jumps
increases and their size decreases, it approximates the so-called
space-time white noise in the sense of convergence of the finite
dimensional distributions. Then, the author proves that the
solutions of the PDEs perturbed by these discrete noises converge in
law (in the sense of finite dimensional distribution convergence)
to the solution of the PDE perturbed by the space-time white noise.

Let us now consider the following one-dimensional quasi-linear stochastic heat equation:
\begin{equation}\label{1nolineal}
\frac{\partial U}{\partial t}(t,x)-\frac{\partial ^2 U}{\partial
x^2} (t,x)=b(U(t,x))+\dot W(t,x),\quad (t,x)\in [0,T]\times [0,1],
 \end{equation}
where $T>0$ stands for a fixed time horizon, $b:\R\rightarrow \R$ is
a globally Lipschitz function and $\dot W$ is the formal notation
for the space-time white noise. We impose some initial condition and
boundary conditions of Dirichlet type, that is:
$$U(0,x)=u_0(x),\; x\in [0,1],$$
$$U(t,0)=U(t,1)=0,\; t\in [0,T],$$
where $u_0:[0,1]\rightarrow \R$ is a continuous function. The random field solution to Equation (\ref{1nolineal}) will be
denoted by $U=\{U(t,x),\; (t,x)\in [0,T]\times[0,1]\}$ and it is interpreted in the {\it{mild}} sense.
More precisely, let $\{W(t,x), \; (t,x)\in [0,T]\times [0,1]\}$ denote a
Brownian sheet on $[0,T]\times [0,1]$, which we suppose to be
defined in some probability space $(\Om,\tf,P)$. For $0\leq t\leq
T$, let $\mathcal{F}_t$ be the $\sigma$-field generated by the
random variables $\{W(s,x), \; (s,x)\in [0,t]\times [0,1]\}$, which
can be conveniently completed, so that the
resulting filtration $\{\mathcal{F}_t,\; t\geq 0\}$ satisfies the usual conditions. Then, a process $U$ is a solution of
(\ref{1nolineal}) if it is $\tf_t$-adapted and the following stochastic integral equation is satisfied:
\begin{align}
U(t,x) = & \int_0^1\! G_t(x,y) u_0(y) dy+\int_0^t\!\!\int_0^1\! G_{t-s}(x,y)\,b(U(s,y)) dy ds \nonumber \\
& \qquad +\int_0^t\!\!\int_0^1\!G_{t-s}(x,y)W(ds,dy), \quad a.s.
\label{mild}
\end{align}
for all $(t,x)\in (0,T]\times (0,1)$, where $G$ denotes the Green function associated to the heat equation in $[0,1]$ with Dirichlet boundary conditions.
We should mention that the stochastic integral in the right-hand side of Equation (\ref{mild}) is a Wiener integral,
which can be understood either in the sense of Walsh \cite{walsh} or in the framework of Da Prato and Zabczyk \cite{dz}.
Besides, existence, uniqueness and pathwise continuity of the solution of (\ref{mild}) are a consequence of \cite[Theorem 3.5]{walsh}.

The aim of our work is to prove that the mild solution of (\ref{1nolineal}) --which is given by the solution of (\ref{mild})-- can be approximated in law,
in the space $\mathcal C([0,T]\times[0,1])$ of continuous functions, by the solution of
\begin{equation}\label{2nolineal}
\frac{\partial U_n}{\partial t}(t,x)-\frac{\partial ^2 U_n}{\partial
x^2} (t,x)=b(U_n(t,x))+\theta_n(t,x),\qquad(t,x)\in [0,T]\times
[0,1],
\end{equation}
with initial condition $u_0$ and Dirichlet boundary conditions,
where $n\in \N$. In this equation, $\theta_n$ will be a noisy input that
approximates the white noise $\dot W$ in the following sense:
\begin{hypothesis}\label{hyp1}
The finite dimensional distributions of the processes
$$\zeta_n(t,x)=\int_0^t\int_0^x \theta_n(s,y) dyds,\quad (t,x)\in [0,T]\times [0,1],$$
converge in law to those of the Brownian sheet 
\end{hypothesis}

Observe that, if the processes $\theta_n$ have square integrable
paths, then the mild form of Equation (\ref{2nolineal}) is given by:
\begin{align}
U_n(t,x)= & \int_0^1\! G_t(x,y) u_0(y) dy+\int_0^t\!\!\int_0^1\!G_{t-s}(x,y)\,b(U_n(s,y)) dy ds \nonumber \\
& \qquad +\int_0^t\!\!\int_0^1\!G_{t-s}(x,y)\theta_n(s,y) dy ds.
\label{mild-approx}
\end{align}
Standard arguments yield existence and uniqueness of solution for Equation
(\ref{mild-approx}) and, furthermore, as it will be detailed later on (see Section \ref{general}),
the solution $U_n$ has continuous trajectories a.s.

In order to state the main result of the paper, let us consider the following hypothesis which, as it will be made explicit in the sequel, will play an essential role:
\begin{hypothesis}\label{hyp2}
For some $q\in [2,3)$, there exists a positive constant $C$ such that, for any $f\in L^q([0,T]\times[0,1])$, it holds:
$$ E\left( \int_0^T \int_0^1 f(t,x) \theta_n(t,x) \, dx dt\right)^2\leq C_q \left( \int_0^T
\int_0^1 |f(t,x)|^{q}\, dx dt\right)^{\frac 2{ q}}.$$
\end{hypothesis}
\begin{hypothesis}\label{hyp3}
There exist $m> 8$ and a positive constant $C$ such that the following is satisfied: for all $s_0,\,s_0'\in [0,T]$ and $x_0,\,x_0'\in
[0,1]$ satisfying $0<s_0<s_0'<2s_0$ and $0<x_0<x_0'<2x_0$, and for any $f\in L^2([0,T]\times [0,1])$, it holds:
$$\sup_{n\geq 1}E\left|\int_{s_0}^{s_0'}\int_{x_0}^{x_0'}f(s,y)\,\theta_n(s,y)dyds\right|^m
\leq C \left( \int_{s_0}^{s_0'}\int_{x_0}^{x_0'}  f(s,y)^2 \, dy ds\right)^{\frac m2}.$$
\end{hypothesis}

We remark that, in Hypothesis \ref{hyp2}, the restriction on the parameter $q$ will be due to the integrability properties of the Green function $G$. On the other hand, in the condition $s_0'<2 s_0$ (resp. $x_0'<2 x_0$) of  Hypothesis \ref{hyp3}, the number $2$ could be replaced by 
any $k>1$. We are now in position to state our main result: 

\begin{teo}\label{quasilinear}
Let $\{\theta_n(t,x),\; (t,x)\in [0,T]\times[0,1]\}$, $n\in\mathbb N$, be
a family of stochastic processes such that $\theta_n\in L^2([0,T]\times [0,1])$ a.s.,  
and such that Hypothesis \ref{hyp1}, \ref{hyp2} and \ref{hyp3} are satisfied. Moreover, assume that $u_0:[0,1]\rightarrow \R$ is continuous and $b:\R\rightarrow \R$ is Lipschitz. 

Then, the family of stochastic processes
$\{U_n,\; n\geq 1\}$ defined as the mild solutions of Equation
(\ref{2nolineal}) converges in law, in the space $\C([0,T]\times
[0,1])$, to the mild solution $U$ of Equation (\ref{1nolineal}).
\end{teo}

Let us point out that, as we will see in Section \ref{general},
Theorem \ref{quasilinear} will be almost an immediate consequence of
the analogous result when taking null initial condition and
nonlinear term (see Theorem \ref{quasilinear2}). Thus, the essential
part of the paper will be concerned to prove the convergence in law,
in the space $\C([0,T]\times [0,1])$, of the solution of
\begin{equation}
\frac{\partial X_n}{\partial t}(t,x)-\frac{\partial ^2 X_n}{\partial
x^2} (t,x)=\theta_n(t,x),\quad (t,x)\in [0,T]\times [0,1],
\label{111}
\end{equation}
with vanishing initial data and Dirichlet boundary conditions, towards the solution of
\begin{equation}
\frac{\partial X}{\partial t}(t,x)-\frac{\partial ^2 X}{\partial
x^2} (t,x)=\dot W(t,x),\quad (t,x)\in [0,T]\times [0,1].
\label{1}
\end{equation}
Observe that the mild solution of
Equations (\ref{111}) and (\ref{1}) can be explicitly written as,
respectively,
\begin{equation}\label{6}
X_n(t,x) = \int_0^t\int_0^1 G_{t-s}(x,y)\, \theta_n(s,y)\; dyds,\;
(t,x)\in [0,T]\times [0,1],
\end{equation}
and
\begin{equation}\label{3}
X(t,x) =  \int_0^t\int_0^1 G_{t-s}(x,y)
W(ds,dy),\; (t,x)\in [0,T]\times [0,1],
\end{equation}
where the latter defines a centered Gaussian process.

\medskip

An important part of the work is also devoted to check that two
interesting particular families of noises verify the hypotheses of
Theorem \ref{quasilinear}. More precisely, consider the following processes:

\begin{enumerate}
\item The \emph{Kac-Stroock processes} on the plane:
\beq
\theta_n(t,x)= n \sqrt{tx} \,(-1)^{N_n(t,x)},
\label{7}
\eeq
where $N_n(t,x):= N(\sqrt{n}t,\sqrt{n}x)$, and $\{N(t,x),\; (t,x)\in [0,T]\times [0,1]\}$
is a standard Poisson process in the plane.
\item The \emph{Donsker kernels}:  Let $\{Z_k,\,k\in\N^2\}$ be an independent family of identically
distributed and centered random variables, with $E(Z_k^2)=1$ for all
$k\in\N^2$, and such that $E(|Z_k|^m)<+\infty$ for all $k\in\N^2$
and some sufficiently large $m\in\N$. For any $n\in\N$, we define the
kernels \beq \theta_n(t,x)=n\sum_{k=(k^1,k^2)\in\N^2}Z_k\cdot {\bf
1}_{[k^1-1,k^1)\times[k^2-1,k^2)}(t n,x n), \;
(t,x)\in[0,T]\times[0,1]. \label{8} \eeq
\end{enumerate}

In the case where $\theta_n$ are the Kac-Stroock processes, it has been proved in \cite{bj2000} that the family of processes
$$\zeta_n(t,x)=\int_0^t\int_0^x \theta_n(s,y) dsdy,\quad n\in \N,$$
converge in law, in the space of continuous functions $\C([0,1]^2)$, to the Brownian sheet.
This result has been inspired by its one-dimensional
counterpart, which is due to Stroock \cite{stroock} and states that the family of
processes
$$Y_\ep(t)=\frac{1}{\ep} \int_0^t (-1)^{N(\frac{s}{\ep^2})} ds,\; t\in [0,1],\; \ep>0,$$
where $N$ stands for a standard Poisson process, converges in law in
$\C([0,1])$ , as $\ep$ tends to $0$, to the standard Brownian motion. Moreover,
it is worth mentioning that Kac (see \cite{kac}) already considered this
kind of processes in order to write the solution of the telegrapher's equation in terms
of a Poisson process.

On the other hand, when $\theta_n$ are the Donsker kernels,
the convergence in law, in the space of continuous
functions, of the processes
$$\zeta_n(t,x)=\int_0^t\int_0^x \theta_n(s,y) dsdy, \quad n\in \N,$$
to the Brownian sheet is a consequence of the extension of Donsker's theorem to the plane (see, for
instance, \cite{wichura}).

\medskip

We should mention at this point that the motivation behind our
results has also been considered by Manthey in \cite{manthey1} and
\cite{manthey2}. Indeed, in the former paper, the author considers
Equation (\ref{111}) with a family of correlated noises
$\{\theta_n,\, n\in \mathbb N \}$ whose integral processes
$$\int_0^t\int_0^x\theta_n(s,y)\,dyds,$$
converge in law (in the sense of finite dimensional distribution
convergence) to the Brownian sheet. Then, sufficient conditions on
the noise processes are specified under which the solution $X_n$ of
(\ref{111}) converges in law, in the sense of the finite dimensional
distribution convergence, to the solution of (\ref{1}). Moreover,
it has also been proved that, whenever the noisy processes are
Gaussian, the convergence in law holds in the space of continuous
functions too; these results have been extended to the quasi-linear
equation (\ref{2nolineal}) in \cite{manthey2}. In this sense, let us
mention that, in an Appendix and for the sake of completeness, we
have added a brief explanation of Manthey's method and showed that
his results do not apply to the examples of noisy inputs that we are
considering in the paper.

Let us also remark that recently there has been an
increasing interest in the study of weak approximation for several
classes of SPDEs (see \cite{deb,debp}). In these references, the
methods for obtaining the corresponding approximation sequences are
based on discretisation schemes for the differential operator
driving the equation, and the rate of convergence of the weak approximations is analysed. Hence, this latter framework differs significantly
from the setting that we have described above. On the other hand, we notice that weak convergence for some classes of SPDEs driven by the Donsker kernels have been considered in the literature; namely, a reduced hyperbolic equation on $\mathbb{R}^2_+$ --which is essentially equivalent to a one-dimensional stochastic wave equation-- has been considered in \cite{cf,fn}, while in \cite{tindel}, the author deals with a stochastic elliptic equation with non-linear drift. 
Furthermore, in \cite{tz}, weak convergence of Wong-Zakai approximations for stochastic evolution equations driven by a finite-dimensional Wiener process has been studied.
Eventually, it is worth commenting
that other type of problems concerning SPDEs driven by Poisson-type
noises have been considered e.g. in
\cite{fournier,erika,leon,mueller,saint}.

The paper is organised as follows. In Section \ref{prel}, we will
present some preliminaries on Equation (\ref{1nolineal}), its
linear form (\ref{1}) and some general results on weak
convergence. In Section \ref{general}, we prove the  results of
convergence for equations (\ref{1}) and (\ref{1nolineal}), so that we end up with the proof of Theorem \ref{quasilinear}. The proof of the fact
that  the Kac-Stroock processes satisfy the hypotheses of Theorem \ref{quasilinear} will be carried out in Section
\ref{poisson}, while the analysis in the case of the Donsker kernels
will be performed at Section \ref{donsker}. Finally, we add
an Appendix where we give the proof of Lemma \ref{lema5} and relate
our results with those of Manthey (\cite{manthey1},
\cite{manthey2}).


\section{Preliminaries}
\label{prel}

As it has been explained in the Introduction, we are concerned with the \emph{mild} solution of the formally-written quasi-linear
stochastic heat equation (\ref{1nolineal}). That is, we consider a real-valued stochastic process $\{U(t,x),\; (t,x)\in [0,T]\times [0,1]\}$, which
we assume to be adapted with respect to the natural filtration generated by the Brownian sheet on $[0,T]\times [0,1]$, such that the following integral equation is satisfied (see (\ref{mild})): for all $(t,x)\in [0,T]\times [0,1]$,
\begin{align}
U(t,x) = & \int_0^1\! G_t(x,y) u_0(y) dy+\int_0^t\!\!\int_0^1\! G_{t-s}(x,y)\,b(U(s,y)) dy ds \nonumber \\
& \qquad +\int_0^t\!\!\int_0^1\!G_{t-s}(x,y)W(ds,dy), \quad a.s.,
\label{mild-bis}
\end{align}
where we recall that $G_t(x,y)$,
$(t,x,y)\in \mathbb{R}_+\times (0,1)^2$, denotes the Green function
associated to the heat equation on $[0,1]$ with Dirichlet boundary
conditions. Explicit formulas for $G$ are well-known, namely:
$$G_t(x,y)=\frac{1}{\sqrt{2\pi t}} \sum_{n=-\infty}^{+\infty} \left(
e^{-\frac{(x-y-2n)^2}{4t}}- e^{-\frac{(x+y-2n)^2}{4t}}\right)$$ or
$$G_t(x,y)=2\sum_{n=1}^\infty \sin(n\pi x) \sin(n\pi y) e^{-n^2 \pi^2
t}.$$ Moreover, it holds that
$$
0\leq G_t(x,y)\leq \frac{1}{\sqrt{2\pi t}} e^{-\frac{(x-y)^2}{4t}},\;
t>0,\; x,y\in [0,1].$$

We have already commented in the Introduction that, in order to prove Theorem \ref{quasilinear},
we will restrict our analysis to the linear version of Equation (\ref{1nolineal}), which is given by
(\ref{1}). Hence, let us consider for the moment  $X=\{X(t,x),\; (t,x)\in
[0,T]\times [0,1]\}$ to be the mild solution of Equation (\ref{1}) with vanishing initial conditions and Dirichlet boundary conditions. This can be 
explicitly written as (\ref{3}). 
Notice that, for any $(t,x)\in (0,T]\times(0,1)$, $X(t,x)$ defines a centered Gaussian random
variable with variance
$$E(X(t,x)^2)=\int_0^t \int_0^1 G_{t-s}(x,y)^2 dy ds.$$
Indeed, by (iii) in Lemma \ref{lema3} below, it holds that
$E(X(t,x)^2)\leq C t^{\frac 1 2}$, where the constant $C>0$ does not
depend on $x$.

\medskip

In the sequel, we will make use of the following result, which is a
quotation of \cite[Lemma B.1]{bms}:
\begin{lema} \label{lema3}
\begin{itemize}
  \item[(i)] Let $\al\in (\frac 3 2,3)$. Then, for all $t\in [0,T]$
  and $x,y\in [0,1]$,
  $$ \int_0^t \int_0^1 |G_{t-s}(x,z)-G_{t-s}(y,z)|^\al dz ds\leq C
  |x-y|^{3-\al}.
  $$
  \item[(ii)] Let $\al\in (1,3)$. Then, for all $s,t\in [0,T]$ such
  that $s\leq t$ and $x\in [0,1]$,
  $$
  \int_0^s \int_0^1 |G_{t-r}(x,y)-G_{s-r}(x,y)|^\al dy dr\leq C
  (t-s)^{\frac{3-\al}{2}}.
  $$
  \item[(iii)] Under the same hypothesis as (ii),
  $$
  \int_s^t \int_0^1 |G_{t-r}(x,y)|^\al dy dr \leq C
  (t-s)^{\frac{3-\al}{2}}.
  $$
\end{itemize}
\end{lema}

\bigskip

Let us recall that we aim to prove that the process $X$ can be approximated
in law, in the space $\C([0,T]\times [0,1])$, by the
family of stochastic processes
\beq
X_n(t,x)=\int_0^t \int_0^1 G_{t-s}(x,y) \theta_n(s,y)\; dy ds, \; (t,x)\in [0,T]\times [0,1],\; n\geq 1,
\label{89}
\eeq
where the processes $\theta_n$ satisfy certain conditions.

In order to prove this convergence in law, we will make use of the following two general results. The first one (Theorem \ref{bill2}) is a tightness criterium on
the plane that generalizes a well-known theorem of Billingsley; it can be found in \cite[Proposition 2.3]{yor}, where it is proved that the hypotheses
considered in the result are stronger than those of the commonly-used criterium of Centsov \cite{centsov}. The second one (Lemma \ref{lema5}) will be used to
prove the convergence of the finite dimensional distributions of $X_n$; though it can be found around in the literature, we have not
been able to find an explicit proof, so that, for the sake of
completeness, we will sketch it in the Appendix.

\begin{teo}\label{bill2}
Let $\{X_n,\,n\in\mathbb N\} $ be a family of random variables
taking values in $\mathcal C([0,T]\times [0,1])$. The family of the
laws of $\{X_n,\,n\in\mathbb N\} $ is tight if  there exist $p',
p>0$, $\delta>2$ and a constant $C$ such that
$$\sup_{n\ge 1 } E|X_n(0,0)|^{p'}<\infty$$
and, for every $t,t'\in[0,T]$ and $x,x'\in[0,1]$,
$$\sup_{n\geq 1}E\left|X_n(t',x')-X_n(t,x)\right|^p\leq
C \left(|x'-x|+|t'-t|\right)^{\delta}.$$
\end{teo}

\begin{lema}\label{lema5}
Let $(F,\|\cdot\|)$ be a normed space and $\{J^n,\; n\in\mathbb N\}$
and $J$  linear maps defined on $F$ and taking values in the space
$L^0(\Om)$ of almost surely finite random variables. Assume that
there exists a positive constant $C$ such that, for any $f\in F$,
\beq \sup_{n\geq 1} E|J^n(f)|\leq C
\|f\|\quad\text{and}\label{30}\eeq \beq
 E|J(f)|\leq C \|f\|,\label{30'}\eeq and that, for
some dense subspace $D$ of $F$, it holds that $J^n(f)$ converges in
law to $J(f)$, as $n$ tends to infinity, for all $f\in D$.

Then, the sequence of random variables $\{J^n(f),\, n\in \N\}$ converges in law
to $J(f)$,  for any $f\in F$.
\end{lema}

\medskip

Eventually, for any real function $X$ defined on $\R^2_+$, and $(t,x),(t',x')\in
\R^2_+$ such that $t\leq t'$ and $x\leq x'$, we will use the
notation $\Del_{t,x}X(t',x')$ for the increment of $X$ over the
rectangle $(t,t']\times (x,x']$:
$$\Del_{t,x}X(t',x')= X(t',x')-X(t,x')-X(t',x)+X(t,x).$$

\section{Proof of the general result}\label{general}

This section is devoted to prove Theorem \ref{quasilinear}. For this, as we have already mentioned, it is convenient to consider, first, the 
linear equation (\ref{1}) together with its mild solution (\ref{3}).

The first step consists in establishing sufficient conditions for a family of processes
$\{\theta_n,\,n\in\mathbb N\}$ in order that the approximation processes $X_n$ (see (\ref{89})) converge, in the sense of finite dimensional distributions,  
to $X$, the solution of (\ref{3}):
\beq 
X(t,x)=\int_0^t\int_0^1 G_{t-s}(x,y)  W(ds,dy).
\label{87} 
\eeq

\begin{prop}\label{dimfin}
Let $\{\theta_n(t,x),\; (t,x)\in [0,T]\times[0,1]\}$, $n\in\mathbb N$, be
a family of stochastic processes such that $\theta_n\in
L^2([0,T]\times [0,1])$ a.s. and such that
Hypothesis \ref{hyp1} and \ref{hyp2} are satisfied.

Then, the finite dimensional distributions
of the processes $X_n$ given by (\ref{89}) converge, as $n$ tends to
infinity, to those of the process defined by (\ref{87}).
\end{prop}

\begin{dem}
We will apply Lemma \ref{lema5} to the following setting: let $q\in
[2,3)$ as in Hypothesis \ref{hyp2} and consider the normed space $(F:=L^q([0,T]\times
[0,1]),\|\cdot\|_q)$, where $\|\cdot\|_q$ denotes the standard norm
in $L^q([0,T]\times [0,1])$. Set
$$J^n(f):=\int_0^T\int_0^1 f(s,y)\theta_n(s,y)\, dy ds,\; \text{ and}$$
$$J(f):=\int_0^T\int_0^1 f(s,y)W( ds, dy),\; f\in F.$$
Then, $J^n$ and $J$ define  linear applications on $F$ and, by
Hypothesis \ref{hyp2}, it holds that
$$\sup_{n\geq 1} E|J^n(f)|\leq C \|f\|_q,$$ for all $f\in
L^q([0,T]\times [0,1])$. The isometry of the Wiener integral
 gives also that
$$E|J(f)|\leq C \|f\|_q,$$for all $f\in
L^q([0,T]\times [0,1])$.
 Moreover, the set $D$ of elementary functions of the form \beq
f(t,x)=\sum_{i=0}^{k-1} f_i \, {\bf 1}_{(t_i,t_{i+1}]}(t){\bf
1}_{(x_i,x_{i+1}]}(x),\label{24}\eeq with $k\geq 1$, $f_i\in \R$,
$0=t_0<t_1<\dots<t_k=T$ and $0=x_0<x_1<\dots<x_k=1$, is dense in
$(F,\|\cdot\|_q)$.

On the other hand, the finite dimensional distributions of $X_n$
converge to those of $X$ if, and only if, for all $m\geq 1$,
$a_1,\dots,a_m\in \R$, $(s_1,y_1),\dots,(s_m,y_m)\in [0,T]\times
[0,1]$, the following convergence in law holds: \beq \sum_{j=1}^m
a_j X_n(s_j,y_j) \bistackrel{\LL}{n\rightarrow
\infty}{\longrightarrow} \sum_{j=1}^m a_j X(s_j,y_j).\label{23}\eeq
This is equivalent to have that $J^n(K)=\int_0^T\int_0^1
K(s,y)\theta_n(s,y)\, dyds$ converges in law, as $n$ tends to
infinity, to $\int_0^T\int_0^1 K(s,y) W(ds,dy)$, where
$$K(s,y):=\sum_{j=1}^m a_j {\bf 1}_{[0,s_j]}(s) G_{s_j-s}(y_j,y).$$
By Lemma \ref{lema3} (iii), the function $K$ belongs to
$L^q([0,T]\times [0,1])$. Hence, owing to
Lemma \ref{lema5}, in order to obtain the convergence (\ref{23}), it
suffices to prove that $J^n(f)$ converges in law to
$J(f)=\int_0^T\int_0^1 f(s,y) W(ds,dy)$, for every elementary
function $f$ of the form (\ref{24}). In fact, if $f$ is such a
function, observe that we have
$$J^n(f)=\sum_{i=0}^{k-1} f_i \int_{t_i}^{t_{i+1}}
\int_{x_i}^{x_{i+1}} \theta_n(s,y)\, dy ds,$$ and this random
variable  converges in law, as $n$ tends to infinity, to
$$\sum_{i=0}^{k-1} f_i \int_{t_i}^{t_{i+1}}
\int_{x_i}^{x_{i+1}} \; W(ds,dy)=\int_0^T\int_0^1 f(s,y) W(ds,dy),$$
because the finite dimensional distributions of $\zeta_n$ converge
to those of the Brownian sheet. \hfill $\Box$

\end{dem}

\medskip

Let us now provide sufficient conditions on $\theta_n$ in order that the
family of laws of the processes $X_n$ is tight in $\mathcal
C([0,T]\times [0,1])$. 

\begin{prop}
Let $\{\theta_n(t,x),\; (t,x)\in [0,T]\times[0,1]\}$, $n\in\mathbb N$, be
a family of stochastic processes such that $\theta_n\in
L^2([0,T]\times [0,1])$ a.s. Suppose that Hypothesis \ref{hyp3} is satisfied. 

Then, the process $X_n$ defined in (\ref{89}) possesses a version with continuous paths and
the family of the laws of $\{X_n,\, n\in \N\}$ is tight in $\mathcal
C([0,T]\times [0,1])$. 
\label{prop2}
\end{prop}

\begin{dem}
It suffices to prove that
\begin{equation}\sup_{n\geq 1}E\left[X_n(t',x')-X_n(t,x)\right]^m\leq C
[|x'-x|^{m\al}+|t'-t|^{\frac{m\al}{2}}],\label{tight} \end{equation}
 for all $\al
\in (0,\frac 12)$,  $t,t'\in [0,T]$ and $x,x'\in [0,1]$. Indeed, if
$m>8$, then it can be found $\al\in (0,\frac12)$ such that
$m\frac{\al}2>2$ and we obtain the existence of a continuous version
of each $X_n$ from Kolmogorov's continuity criterium in the plane. Furthermore, 
by Theorem \ref{bill2}, we also obtain the tightness of the laws of $X_n$ in $\mathcal C([0,T]\times [0,1])$.

Set $H(t,x;s,y):= {\bf 1}_{[0,t]}(s) G_{t-s}(x,y)$. We will need
to estimate the moment of order $m$, for some  $m>8$, of the
quantity
$$X_n(t',x')-X_n(t,x)=\int_0^T\int_0^1 [H(t',x';s,y)-H(t,x;s,y)]
\theta_n(s,y)\; dy ds,$$ for $t,t'\in [0,T]$ and $x,x'\in [0,1]$.
Moreover, the right-hand side of the above equality can be written
in the form $\Del_{0,0} Y_n(T,1)$, where the process $Y_n$, which
indeed depends on $t,t',x,x'$, is defined by
$$Y_n(s_0,x_0):=\int_0^{s_0}\int_0^{x_0} [H(t',x';s,y)-H(t,x;s,y)]
\theta_n(s,y)\; dy ds,\; (s_0,x_0)\in [0,T]\times [0,1].$$ Hence,
inequality (\ref{tight}) is equivalent to prove that
$$E(\Del_{0,0} Y_n(T,1))^m \leq
C [|x'-x|^{m\al}+|t'-t|^{\frac{m\al}{2}}],$$ for all $\al \in
(0,\frac 12)$ and $n\geq 1$. By \cite[Lemma 3.2]{bjt} (in the
statement of this lemma, it is supposed that $m$ is an even integer
number, but this assumption is not used in its proof), it suffices
to prove that there exist $\gam>0$ and $C>0$ such that, for all
$s_0,s_0'\in [0,T]$ and $x_0,x_0'\in [0,1]$ satisfying
$0<s_0<s_0'<2s_0$ and $0<x_0<x_0'<2x_0$, then \beq \sup_{n\geq
1}E(\Del_{s_0,x_0}Y_n(s_0',x_0'))^m \leq C \left[
|t'-t|^{m\al}+|x'-x|^{\frac{m\al}{2}} \right] (s_0'-s_0)^{m\gam}
(x_0'-x_0)^{m\gam}.\label{18}\eeq
%
%
%
%
By Hypothesis \ref{hyp3} for the particular case of $f(s,y)=H(t',x';s,y)
-H(t,x;s,y)$, we obtain
\begin{align*}
& \sup_{n\geq 1}E(\Del_{s_0,x_0}Y_n(s_0',x_0'))^m \\
&\quad \quad  \leq C \left( \int_0^T\int_0^1 {\bf
1}_{[s_0,s_0']}(s){\bf 1}_{[x_0,x_0']}(y)
|H(t',x';s,y)-H(t,x;s,y)|^2 \, dy ds\right)^{\frac m2}.
\end{align*}
Let $p\in (1,\frac 32)$ and $q>1$ such that $\frac 1p + \frac 1q
=1$. Then, by H\"{o}lder's inequality and the definition of $H$,
\begin{align}
& \sup_{n\geq 1}E(\Del_{s_0,x_0}Y_n(s_0',x_0'))^m  \nonumber \\
& \leq C \left( \int_0^T\int_0^1 {\bf 1}_{[s_0,s_0']}(s){\bf
1}_{[x_0,x_0']}(y) \, dy ds \right)^{\frac{m}{2q}} \left(
\int_0^T\int_0^1 |H(t',x';s,y)-H(t,x;s,y)|^{2p} \, dy
ds\right)^{\frac{m}{2p}}\nonumber\\
&  \leq C (x_0'-x_0)^{\frac{m}{2q}}
(s_0'-s_0)^{\frac{m}{2q}} \nonumber \\
& \quad \quad \times \left( \int_0^t\int_0^1
|G_{t'-s}(x',y)-G_{t-s}(x,y)|^{2p} \, dy ds + \int_t^{t'} \int_0^1
|G_{t'-s}(x',y)|^{2p}\, dy ds \right)^{\frac{m}{2p}}. \label{22}
\end{align}
By Lemma \ref{lema3}, the last term in the right-hand side of
(\ref{22}) can be bounded, up to some constant, by
$$\left( |x-x'|^{3-2p}+|t-t'|^{\frac{3-2p}{2}}\right)^{\frac{m}{2p}}
\leq C \left(
|x-x'|^{\frac{m(3-2p)}{2p}}+|t-t'|^{\frac{m(3-2p)}{4p}}\right).$$
Therefore, if we plug this bound in (\ref{22}) and we take
$\al=\frac{3-2p}{2p}$ and $\gam=\frac{1}{2q}$, then we have proved
(\ref{18}), because $p\in (1,\frac 32)$ is arbitrary. \hfill $\Box$
\end{dem}

\begin{obs}
As it can be deduced from the first part of the proof of Proposition \ref{prop2}, the restriction $m>8$ has to be considered 
in order to be able to apply Theorem \ref{bill2} and Kolmogorov's continuity criterium.
\end{obs}

As a consequence of Propositions \ref{dimfin} and \ref{prop2}, we can state the following result
on convergence in law for the processes $X_n$:

\begin{teo}\label{linear}
Let $\{\theta_n(t,x),\; (t,x)\in [0,T]\times[0,1]\}$, $n\in\mathbb N$, be
a family of stochastic processes such that $\theta_n\in
L^2([0,T]\times [0,1])$ a.s. Assume that Hypothesis \ref{hyp1}, \ref{hyp2} and \ref{hyp3} are satisfied.

Then, the family of stochastic processes $\{X_n,\; n\geq 1\}$ defined in (\ref{89}) converges in law, as $n$
tends to infinity in the space $\C([0,T]\times [0,1])$, to the Gaussian process $X$ given by (\ref{87}).
\end{teo}

\medskip

We can eventually extend the above result to the quasi-linear Equation (\ref{1nolineal}), so that we end up with the proof 
of Theorem \ref{quasilinear}. This will be an immediate consequence of the above theorem and the next general result:

\begin{teo}\label{quasilinear2}
Let $\{\theta_n(t,x),\; (t,x)\in [0,T]\times[0,1]\}$, $n\in\mathbb N$, be
a family of stochastic processes such that $\theta_n\in
L^2([0,T]\times [0,1])$ a.s. Assume that $u_0:[0,1]\rightarrow \R$ is a continuous function and $b:\R\rightarrow \R$ is Lipschitz. Moreover, 
suppose that the family of stochastic processes $\{X_n,\; n\geq 1\}$ defined in (\ref{89}) converges in law, as $n$ tends to infinity in the space
$\C([0,T]\times [0,1])$, to the Gaussian process $X$ given by (\ref{87}). 

Then, the family of stochastic processes $\{U_n,\; n\geq
1\}$ defined as the mild solutions of Equation (\ref{2nolineal})
converges in law, in the space $\C([0,T]\times [0,1])$, to the mild
solution $U$ of Equation (\ref{1nolineal}).
\end{teo}

\begin{dem}
Let us first recall that we denote by $U=\{U(t,x),\; (t,x)\in
[0,T]\times [0,1]\}$ the unique mild solution of Equation
(\ref{1nolineal}), which means that $U$ fulfils
\begin{align*}
U(t,x) = & \int_0^1\! G_t(x,y) u_0(y) dy+\int_0^t\!\!\int_0^1\! G_{t-s}(x,y)\,b(U(s,y)) dy ds  \\
& \qquad +\int_0^t\!\!\int_0^1\!G_{t-s}(x,y)W(ds,dy), \quad a.s.
\end{align*}
The approximation sequence is denoted by $\{U_n,\; n\in \N\}$, where
$U_n=\{U_n(t,x),\; (t,x)\in [0,T]\times [0,1]\}$ is a stochastic
process satisfying
\begin{align*}
U_n(t,x) = & \int_0^1\! G_t(x,y) u_0(y) dy+\int_0^t\!\!\int_0^1\! G_{t-s}(x,y)\,b(U_n(s,y)) dy ds  \\
& \qquad +\int_0^t\!\!\int_0^1\!G_{t-s}(x,y) \theta_n(s,y) dy ds,
\quad a.s. 
\end{align*}
where the noisy input $\theta_n$ has square integrable paths, a.s.

\medskip


%


Using the properties of the Green function (see Lemma \ref{lema3}),
the fact that $\theta_n\in L^2([0,T]\times [0,1])$ a.s., together
with a Gronwall-type  argument, we obtain that $U_n$ has continuous
paths a.s., for all $n\in \N$.

Next, for each continuous function $\eta:[0,T]\times
[0,1]\longrightarrow \mathbb R$, consider the following
(deterministic) integral equation:
$$z_{\eta}(t,s)=\int_0^1\! G_t(x,y) u_0(y) dy+\int_0^t\!\!\int_0^1\!
G_{t-s}(x,y)\,b(z_{\eta}(s,y)) dy ds+\eta(t,x).$$ 
As before, by the properties of $G$ and the assumptions on $u_0$ and $b$, it can be checked
that this equation possesses a unique continuous solution.

 Now, we
will prove that the map
$$\psi:\mathcal C([0,T]\times [0,1])\longrightarrow \mathcal
C([0,T]\times[0,1])$$
$$\phantom{xx}\eta\phantom{xxx}\longrightarrow \phantom{xx}z_{\eta}$$
is continuous with respect to the usual topology on this space.
Indeed, given $\eta_1,\,\,\eta_2\in \mathcal C([0,T]\times[0,1])$,
we have that
\begin{align}
& |z_{_{\eta_1}}(t,x)-z_{_{\eta_2}}(t,x)| \nonumber \\
&\qquad \le \int_0^t\!\!\int_0^1
\!G_{t-s}(x,y) \,\left|b(z_{_{\eta_1}}(s,y))-b(z_{_{\eta_2}}(s,y))\right|dyds+|\eta_1(t,x)-\eta_2(t,x)|\nonumber \\
&\qquad \le L\int_0^t\!\!\int_0^1 \!G_{t-s}(x,y)
\,\left|z_{_{\eta_1}}(s,y)-z_{_{\eta_2}}(s,y)\right|dyds+|\eta_1(t,x)-\eta_2(t,x)|,
\label{func}
\end{align}
where $L$ is the Lipschitz constant of the function $b$.

For a given  $f\in \mathcal C([0,T]\times [0,1])$, we introduce the
following norms:
$$\|f\|_t=\max_{s\in[0,\,t],\,x\in[0,\,1]}|f(s,x)|.$$
By using this notation, we deduce that inequality (\ref{func})
implies that, for any $t\in[0,T]$,
$$\|z_{_{\eta_1}}-z_{_{\eta_2}}\|_t\le L\int_0^t\overline{
G}(t-s)\,\,\|z_{_{\eta_1}}-z_{_{\eta_2}}\|_s\, ds
+\|\eta_1-\eta_2\|_{_T},$$ where
$$\overline{G}(s):=\sup_{x\in[0,\,1]}\int_0^1 G_s(x,y) dy\le \sup_{x\in[0,\,1]}\int_0^1
\frac1{\sqrt{2\pi s}}e^{-\frac{(x-y)^2}{4s}} dy\le C.$$
Applying now Gronwall's lemma, we obtain that there exists a finite
constant $A>0$ such that
$$\|z_{_{\eta_1}}-z_{_{\eta_2}}\|_{_T}\le A \,\|\eta_1-\eta_2\|_{_T},$$
and, therefore, the map $\psi$ is continuous.

Consider now
 $$X_n(t,x)=\int_0^t\int_0^1 G_{t-s}(x,y)\theta_n(s,y) dy ds$$
 and
$$X(t,x)=\int_0^t\int_0^1 G_{t-s}(x,y)W(ds,dy).$$
By hypothesis, we have that $X_n$ converges in law in $\mathcal
C([0,T]\times[0,1])$ to $X$, as $n$ goes to infinity. On the other
hand, we have
$$U_n=\psi(X_n)\quad \text{and}\quad U=\psi(X),$$
and hence the continuity of $\psi$ implies the convergence in law of
$U_n$ to $U$ in $\mathcal C([0,T]\times[0,1])$. \hfill $\Box$
\end{dem}
\section{Convergence in law for the Kac-Stroock processes}
\label{poisson}

This section is devoted to prove that the hypotheses of Theorem
\ref{quasilinear} are satisfied in the case where
the approximation family is defined in terms of the Kac-Stroock
process $\theta_n$ set up in (\ref{7}). That is, 
\beq X_n(t,x)= n
\int_0^t \int_0^1 G_{t-s}(x,y) \sqrt{sy} (-1)^{N_n(s,y)}\; dy ds.
\label{17} 
\eeq 
First, we notice that Hypothesis \ref{hyp1} has been proved in \cite{bj2000}. 

The following proposition states that Hypothesis \ref{hyp2} is satisfied in this particular situation.
\begin{prop}
Let $\theta_n$ be the Kac-Strock processes. Then, for all $p>1$,
there exists a positive constant $C_p$ such that \beq E\left(
\int_0^T \int_0^1 f(t,x) \theta_n(t,x) \, dx dt\right)^2\leq C_p
\left( \int_0^T \int_0^1 |f(t,x)|^{2p}\, dx dt\right)^{\frac 1 p},
\label{25} \eeq for any $f\in L^{2p}([0,T]\times [0,1])$ and all
$n\geq 1$. \label{prop1}
\end{prop}
The proof of this proposition is based on the following technical
lemma:
\begin{lema}
Let $f\in L^2([0,T]\times [0,1])$ and $\al\geq 1$. Then, for any
$u,u'\in (0,1)$ satisfying that $0<u<u'\leq 2^\al u$,
$$E\left(  \int_0^T \int_u^{u'} f(t,x) \theta_n(t,x)\, dx dt\right)^2
\leq \frac 34 \left( 2^{\al+1} -1 \right) \int_0^T \int_u^{u'}
f^2(t,x) \, dx dt,$$ for all $n\geq 1$. \label{lema4}
\end{lema}

\begin{dem}
First, we observe that
\begin{align}
E\left( \int_0^T \int_u^{u'} f(t,x) \theta_n(t,x)\, dx dt\right)^2 =
&  2 n^2 \int_0^T \int_u^{u'}\int_0^T \int_u^{u'} f(t_1,x_1)
f(t_2,x_2) \sqrt{t_1t_2x_1x_2} \nonumber \\
& \quad \times  E\left[ (-1)^{N_n(t_1,x_1)+N_n(t_2,x_2)}\right] {\bf
1}_{\{ t_1\leq t_2\}} dx_2 dt_2dx_1dt_1. 
\label{9} 
\end{align} 
The
expectation appearing in (\ref{9}) can be computed as it has been
done in the proof of \cite[Lemma 3.1]{bjt} (see also \cite[Lemma
3.2]{bj2000}). More precisely, one writes the sum
$N_n(t_1,x_1)+N_n(t_2,x_2)$ as a suitable sum of rectangular
increments of $N_n$ and applies that, if $Z$ has a Poisson
distribution with parameter $\lam$, then $E\left[(-1)^Z\right]
=\exp(-2\lam)$. Hence, the term in the right-hand side of (\ref{9})
admits a decomposition of the form $I_1+I_2$, where
\begin{align*}
I_1 = &  2 n^2 \int_0^T \int_u^{u'}\int_0^T \int_u^{u'} f(t_1,x_1)
f(t_2,x_2) \sqrt{t_1t_2x_1x_2}  \\
& \quad \quad \quad \times \exp\left\{-2n[(t_2-t_1)x_2 +
(x_2-x_1)t_1]\right\} {\bf 1}_{\{ t_1\leq t_2\}} {\bf 1}_{\{ x_1\leq
x_2\}} dx_2 dt_2dx_1dt_1, 
\end{align*}
\begin{align*}
I_2 = &  2 n^2 \int_0^T \int_u^{u'}\int_0^T \int_u^{u'} f(t_1,x_1)
f(t_2,x_2) \sqrt{t_1t_2x_1x_2}  \\
& \quad \quad \quad \times \exp\left\{-2n[(t_2-t_1)x_2 +
(x_1-x_2)t_1]\right\} {\bf 1}_{\{ t_1\leq t_2\}} {\bf 1}_{\{ x_2\leq
x_1\}} dx_2 dt_2dx_1dt_1. 
\end{align*} 
Let us apply the
inequality $ab\leq \frac 1 2 (a^2+b^2)$, $a,b\in \R$, so that we
have $I_1\leq I_{11}+I_{12}$, where the latter terms are defined by
\begin{align*}
I_{11} = &   n^2 \int_0^T \int_u^{u'}\int_0^T \int_u^{u'}
f^2(t_1,x_1) \,
t_1 x_1  \\
& \quad \quad \quad \times \exp\left\{-2n[(t_2-t_1)x_2 +
(x_2-x_1)t_1]\right\} {\bf 1}_{\{ t_1\leq t_2\}} {\bf 1}_{\{ x_1\leq
x_2\}} dx_2 dt_2dx_1dt_1, \end{align*}
\begin{align*}
I_{12} = &   n^2 \int_0^T \int_u^{u'}\int_0^T \int_u^{u'}
f^2(t_2,x_2) \,
t_2 x_2  \\
& \quad \quad \quad \times \exp\left\{-2n[(t_2-t_1)x_2 +
(x_2-x_1)t_1]\right\} {\bf 1}_{\{ t_1\leq t_2\}} {\bf 1}_{\{ x_1\leq
x_2\}} dx_2 dt_2dx_1dt_1. \end{align*} In order to deal with the
term $I_{11}$, we will use the fact that $\exp\{-2n
(t_2-t_1)x_2\}\leq \exp\{-2n (t_2-t_1)x_1\}$, for $x_1\leq x_2$, and
then integrate with respect to $t_2,x_2$. Thus
\begin{align}
I_{11} & \leq    n^2 \int_0^T \int_u^{u'}\int_0^T \int_u^{u'}
f^2(t_1,x_1) \,
t_1 x_1 \nonumber \\
& \quad \quad \quad \times \exp\left\{-2n[(t_2-t_1)x_1 +
(x_2-x_1)t_1]\right\} {\bf 1}_{\{ t_1\leq t_2\}} {\bf 1}_{\{ x_1\leq
x_2\}} dx_2
dt_2dx_1dt_1 \nonumber \\
& \leq \frac 1 4  \int_0^T \int_u^{u'} f^2(t_1,x_1) \, dx_1 dt_1.
\label{12}
\end{align} Concerning the term $I_{12}$, we use similar arguments
as before and, moreover, we apply the fact that, for $x_1,x_2\in
[u,u')$, then $x_2<2^\al x_1$. Hence
\begin{align}
I_{12} & \leq    n^2 \int_0^T \int_u^{u'}\int_0^T \int_u^{u'}
f^2(t_2,x_2) \,
t_2 x_2 \nonumber \\
& \quad \quad \quad \times \exp\left\{-2n[(t_2-t_1)x_1 +
(x_2-x_1)t_2]\right\} {\bf 1}_{\{ t_1\leq t_2\}} {\bf 1}_{\{ x_1\leq
x_2\}} dx_2
dt_2dx_1dt_1 \nonumber \\
& \leq  2^\al  n^2 \int_0^T \int_u^{u'}\int_0^T \int_u^{u'}
f^2(t_2,x_2) \,
t_2 x_1 \nonumber \\
& \quad \quad \quad \times \exp\left\{-2n[(t_2-t_1)x_1 +
(x_2-x_1)t_2]\right\} {\bf 1}_{\{ t_1\leq t_2\}} {\bf 1}_{\{ x_1\leq
x_2\}} dx_2
dt_2dx_1dt_1 \nonumber \\
& \leq 2^{\al-2}  \int_0^T \int_u^{u'} f^2(t_2,x_2) \, dx_2 dt_2.
\label{13}
\end{align}
The analysis of the term $I_2$ is slightly more involved. Namely,
notice first that $I_2\leq I_{21}+I_{22}$, where
\begin{align*}
I_{21} = &   n^2 \int_0^T \int_u^{u'}\int_0^T \int_u^{u'}
f^2(t_1,x_1) \,
t_1 x_1  \\
& \quad \quad \quad \times \exp\left\{-2n[(t_2-t_1)x_2 +
(x_1-x_2)t_1]\right\} {\bf 1}_{\{ t_1\leq t_2\}} {\bf 1}_{\{ x_2\leq
x_1\}} dx_2 dt_2dx_1dt_1, \end{align*}
\begin{align*}
I_{22} = &   n^2 \int_0^T \int_u^{u'}\int_0^T \int_u^{u'}
f^2(t_2,x_2) \,
t_2 x_2  \\
& \quad \quad \quad \times \exp\left\{-2n[(t_2-t_1)x_2 +
(x_1-x_2)t_1]\right\} {\bf 1}_{\{ t_1\leq t_2\}} {\bf 1}_{\{ x_2\leq
x_1\}} dx_2 dt_2dx_1dt_1. \end{align*} For the term $I_{12}$, we
simply use that, by hypothesis, $x_1\leq 2^\al x_2$, and we
integrate with respect to $t_2, x_2$, so that we end up with
\begin{align}
I_{21} & \leq  2^\al  n^2 \int_0^T \int_u^{u'}\int_0^T \int_u^{u'}
f^2(t_1,x_1) \,
t_1 x_2 \nonumber \\
& \quad \quad \quad \times \exp\left\{-2n[(t_2-t_1)x_2 +
(x_1-x_2)t_1]\right\} {\bf 1}_{\{ t_1\leq t_2\}} {\bf 1}_{\{ x_2\leq
x_1\}} dx_2
dt_2dx_1dt_1 \nonumber \\
& \leq 2^{\al-2}  \int_0^T \int_u^{u'} f^2(t_1,x_1) \, dx_1 dt_1.
\label{14}
\end{align}
The term $I_{22}$ is much more delicate. Namely, taking into account
the integration's region in $I_{22}$ as well as the fact that
$x_1-x_2\leq (2^\al-1) x_2$ (because $x_1\leq 2^\al x_2$), it holds
\begin{align*}
2(t_2-t_1)x_2 + 2(x_1-x_2)t_1 & \geq (t_2-t_1)x_2 +
\frac{1}{2^\al-1}(t_2-t_1)(x_1-x_2) + \frac{1}{2^\al-1}(x_1-x_2)t_1
\\
& =(t_2-t_1)x_2 +\frac{1}{2^\al-1}(x_1-x_2)t_2.
\end{align*}
Therefore,
\begin{align}
I_{22} & \leq    n^2 \int_0^T \int_u^{u'}\int_0^T \int_u^{u'}
f^2(t_2,x_2) \,
t_2 x_2 \nonumber \\
& \quad \quad \quad \times \exp\left\{-n[(t_2-t_1)x_2 +
\frac{1}{2^\al-1} (x_1-x_2)t_2]\right\} {\bf 1}_{\{ t_1\leq t_2\}}
{\bf 1}_{\{ x_2\leq x_1\}} dx_2
dt_2dx_1dt_1 \nonumber \\
& \leq (2^\al-1)  \int_0^T \int_u^{u'} f^2(t_2,x_2) \, dx_2 dt_2,
\label{15}
\end{align}
where the latter expression has been obtained after integrating with
respect to $t_1, x_1$.

We conclude the proof by putting together (\ref{12})-(\ref{15}).
\hfill $\Box$
\end{dem}

\bigskip

\textit{Proof of Proposition \ref{prop1}:} Let us consider the
following dyadic-type partition of $(0,1]$:
$$(0,1]=\bigcup_{k=0}^\infty (a_{k+1},a_k],$$
with $a_k=\frac{1}{2^{k\al}}$, for some $\al\geq 1$. In particular,
observe that $a_k-a_{k+1}=\frac{2^\al-1}{2^{(k+1)\al}}$ and we are
in position to apply Lemma \ref{lema4}: for all $k\geq 0$,
$$E \left(\int_0^T\int_{a_{k+1}}^{a_k} f(t,x) \theta_n(t,x)\, dx
dt\right)^2\leq \frac 34 (2^{\al+1}-1) \int_0^T\int_{a_{k+1}}^{a_k}
f(t,x)^2 \, dx dt.$$ Therefore, we have the following estimations:
\begin{align}
E\left( \int_0^T \int_0^1 f(t,x) \theta_n(t,x) \, dx dt\right)^2 & =
E\left( \sum_{k=0}^\infty \int_0^T \int_{a_{k+1}}^{a_k} f(t,x)
\theta_n(t,x) \, dx dt\right)^2 \nonumber \\
& \leq  \sum_{k=0}^\infty 2^{k+1} E\left(\int_0^T
\int_{a_{k+1}}^{a_k}
 f(t,x) \theta_n(t,x) \, dx dt\right)^2 \nonumber \\
& \leq \frac 34 (2^{\al+1}-1) \sum_{k=0}^\infty
2^{k+1}\int_0^T\int_{a_{k+1}}^{a_k} f(t,x)^2 \, dx dt. \label{16}
\end{align}
Let $p,q>1$ be such that $\frac 1 p + \frac 1 q =1$. Then, applying
H\"{o}lder's inequality, the last term of (\ref{16}) can be bounded by
\begin{align}
 & \frac 34 (2^{\al+1}-1) \sum_{k=0}^\infty 2^{k+1} \left( \int_0^T \int_{a_{k+1}}^{a_k}
 |f(t,x)|^{2p} \, dx dt \right)^{\frac 1 p} (a_k-a_{k+1})^{\frac 1
 q}\nonumber \\
& \quad \leq \frac 34 (2^{\al+1}-1) \left( \int_0^T \int_0^1
|f(t,x)|^{2p} \, dx dt \right)^{\frac 1 p}\sum_{k=0}^\infty 2^{k+1}
\frac{(2^\al-1)^{\frac
1 q}}{2^{(k+1)\frac{\al}{q}}}\nonumber \\
& \quad \leq \frac 34 (2^{\al+1}-1)(2^\al-1)^{\frac 1q} \left(
\int_0^T \int_0^1 |f(t,x)|^{2p} \, dx dt \right)^{\frac 1
p}\sum_{k=0}^\infty \frac{1}{2^{(k+1)\left( \frac{\al}{q}-1
\right)}} \label{31}
\end{align}
and this series is convergent whenever we take $\al$ such that
$\al>q$. Hence, expression (\ref{31}) may be bounded by
$$\frac 32 (2^{\al+1}-1) \frac{(2^\al-1)^{\frac 1q}}{2^{\frac \al q -2}} \left(
\int_0^T \int_0^1 |f(t,x)|^{2p} \, dx dt \right)^{\frac 1 p},$$
which implies that the proof is complete.
\hfill $\Box$

\medskip

\begin{obs}
It is worth noticing that, in the statement of Proposition
\ref{prop1}, we have not been able to obtain the validity of the
result for $p=1$. Indeed, as it can be deduced from its proof, the
constant $C_p$ in (\ref{25}) blows up when $p\rightarrow 1$ (because
$q\rightarrow \infty$, so $\al\rightarrow \infty$). \label{obs1}
\end{obs}

By Proposition \ref{dimfin}, a consequence of Proposition \ref{prop1} is that the finite dimensional
distributions of $X_n$ (see (\ref{17})) converge, as $n$ tends to infinity, to those of
$$X(t,x)=  \int_0^t \int_0^1 G_{t-s}(x,y) W(ds,dy).$$

\medskip

In order to prove that Theorem \ref{quasilinear} applies for the Kac-Stroock processes, 
it only remains to verify that Hypothesis \ref{hyp3} is
satisfied. In fact, this is given by the following result:

\begin{prop}
Let $\theta_n$ be the Kac-Stroock kernels. Then, for any even
$m\in\mathbb N$,  there exists a positive constant $C_m$ such that,
for all $s_0,\,s_0'\in [0,T]$ and $x_0,\,x_0'\in [0,1]$ satisfying
$0<s_0<s_0'<2s_0$ and $0<x_0<x_0'<2x_0$, we have that
$$\sup_{n\geq
1}E\left(\int_{s_0}^{s_0'}\int_{x_0}^{x_0'}f(s,y)\,\theta_n(s,y)dyds\right)^m
\leq C_m \left( \int_{s_0}^{s_0'}\int_{x_0}^{x_0'}  f(s,y)^2 \, dy
ds\right)^{\frac m2},
$$
 for any
$f\in L^2([0,T]\times [0,1])$.

\label{prop3}
\end{prop}

\begin{dem}
To begin with, define
$$Z_n(s_0,x_0):=\int_0^{s_0} \int_0^{x_0} f(s,y)\theta_n(s,y)\, dyds$$
and 
observe that we can apply the same arguments as in
the proof of \cite[Lemma 3.3]{bjt} (see p. 324 therein) in order to
obtain the following estimate:
\begin{align*}
E(\Del_{s_0,x_0} Z_n(s_0',x_0'))^m \leq &  m! n^m \int_{[0,T]^m\times
[0,1]^m}  \prod_{i=1}^m \left( {\bf 1}_{[s_0,s_0']}(s_i){\bf
1}_{[x_0,x_0']}(y_i)
f(s_i,y_i) \sqrt{s_i y_i}\right) \nonumber \\
&\quad \quad  \times
\exp\left\{-n[(s_m-s_{m-1})y_{(m-1)}+\cdots+(s_2-s_1)y_{(1)}]\right\}
\nonumber \\
&\quad \quad  \times
\exp\left\{-n[(y_{(m)}-y_{(m-1)})s_{m-1}+\cdots+(y_{(2)}-y_{(1)})s_{1}]\right\}
\nonumber \\
& \quad  \quad \times {\bf 1}_{\{ s_1\leq \dots\leq s_m\}} ds_1
\cdots ds_m dy_1\cdots dy_m,
\end{align*}
where $y_{(1)},\dots,y_{(m)}$ denote the variables $y_1,\dots,y_m$
ordered increasingly. Hence
\begin{align}
E(\Del_{s_0,x_0}Z_n(s_0',x_0'))^m \leq & 2^m (s_0 x_0)^{\frac m2} m!
n^m \int_{[0,T]^m\times [0,1]^m}  \prod_{i=1}^m\left( {\bf
1}_{[s_0,s_0']}(s_i){\bf 1}_{[x_0,x_0']}(y_i)
f(s_i,y_i)\right) \nonumber \\
&\quad \quad  \times \exp\left\{-n
x_0[(s_m-s_{m-1})+\cdots+(s_2-s_1)]\right\}
\nonumber \\
&\quad \quad  \times \exp\left\{-n
s_0[(y_{(m)}-y_{(m-1)})+\cdots+(y_{(2)}-y_{(1)})]\right\}
\nonumber \\
& \quad  \quad \times {\bf 1}_{\{ s_1\leq \dots\leq s_m\}} ds_1
\cdots ds_m dy_1\cdots dy_m. \label{19}
\end{align}
Notice that in (\ref{19}) we have not been able to order the
variables $y_1,\dots,y_m$, because neither the function
$(s,y)\mapsto f(s,y)$ factorizes nor $(y_1\dots,y_m)\mapsto
f(s_1,y_1)\cdots f(s_m,y_m)$ is symmetric. However, the fact that
the variables $s_i$ are ordered determines $\frac m2$ couples
$(s_1,s_2),\,(s_3,s_4)\dots,(s_{m-1},s_m)$, such that the second
element in each couple is greater than or equal to the first one.
Concerning the variables $y_i$, we also have $\frac m2$ couples
$(y_{(1)},y_{(2)}),\dots, (y_{(m-1)},y_{(m)})$ satisfying
the same property.

The key point of the proof relies in factorizing the product in the
first part of the right-hand side of (\ref{19}) into two convenient
products:
$$\prod_{j=1}^{\frac m2}\left(
{\bf 1}_{[s_0,s_0']}(s_{i_j}){\bf 1}_{[x_0,x_0']}(y_{i_j})
f(s_{i_j},y_{i_j})\right) \prod_{k=1}^{\frac m2}\left( {\bf
1}_{[s_0,s_0']}(s_{r_k}){\bf 1}_{[x_0,x_0']}(y_{r_k})
f(s_{r_k},y_{r_k})\right),$$ where $\mathcal{I}=\{i_j,\,
j=1,\dots,\frac m2\}$ and $\mathcal{R}=\{r_k,\, k=1,\dots,\frac
m2\}$ are two disjoint subsequences of $\{1,\dots,m\}$. In
particular, it holds that $\mathcal{I}\uplus
\mathcal{R}=\{1,\dots,m\}$. These subsequences will be chosen using
the following rule: any couple $(s_i,s_{i+1})$ will contain an
element of the form $s_{i_j}$ and one of the form $s_{r_k}$, and any
couple $(y_{(i)},y_{(i+1)})$ will contain an element of the form
$y_{i_j}$ and one of the form $y_{r_k}$. For this, we will split the
$m$ elements $f(s_1,y_1),\dots,f(s_m,y_m)$ in two groups of $\frac
m2$ elements:
$$A=\{f(s_{i_1},y_{i_1}),\dots,f(s_{i_{\frac m2}},y_{i_{\frac
m2}})\},$$
$$B=\{f(s_{r_1},y_{r_1}),\dots,f(s_{r_{\frac m2}},y_{r_{\frac
m2}})\}.$$ In order to determine the elements of each group, and
such that the above condition is satisfied, we proceed by an
iterative method: we will start with an element of $A$ and we will
associate to it an element of $B$ satisfying what we want; then, to
the latter element of $B$ we will associate a suitable element of
$A$, and so on. More precisely, we start, say, with
$f(s_{i_1},y_{i_1})=f(s_1,y_1)$. Then, if at any step of the
iteration procedure we have an element $f(s_{i_j},y_{i_j})\in A$, we
will associate to it an element $f(s_{r_k},y_{r_k})\in B$ in such a
way that $\{s_{i_j},s_{r_k}\}$ forms one of the couples
$(s_i,s_{i+1})$. On the other hand, if at any step of the iteration
procedure we have an element $f(s_{r_k},y_{r_k})\in B$, then we will
associate to it $f(s_{i_j},y_{i_j})\in A$ such that
$\{y_{i_j},y_{r_k}\}$ determines one of the couples
$(y_{(i)},y_{(i+1)})$. The only thing that remains to be clarified
is what we are going to do in case that, at some step, we end up
with and element of $A$ or $B$ which has already appeared before. In
this case, we do not take the latter element, but another one which
has not been chosen by now.

Let us illustrate the above-described procedure by considering a
particular example: let $m=8$ and assume that we fix $y_1,\dots,y_8$
in such a way that
$$y_8<y_5<y_4<y_7<y_1<y_6<y_2<y_3,$$
that is:
$$y_{(1)}=y_8,\, y_{(2)}=y_5,\, y_{(3)}=y_4,\, y_{(4)}=y_7,$$
$$y_{(5)}=y_1,\, y_{(6)}=y_6,\, y_{(7)}=y_2,\, y_{(8)}=y_3.$$
Recall that we assume that $s_1\leq \dots \leq s_8$. We start with
$f(s_1,y_1)\in A$. Then, the iteration sequence will be the
following:
$$f(s_1,y_1)\longrightarrow f(s_2,y_2) \longrightarrow f(s_3,y_3)
\longrightarrow f(s_4,y_4) $$ $$\longrightarrow f(s_7,y_7)
\longrightarrow f(s_8,y_8) \longrightarrow f(s_5,y_5)
\longrightarrow f(s_6,y_6)$$ Thus,
$A=\{f(s_1,y_1),f(s_3,y_3),f(s_7,y_7),f(s_5,y_5)\}$ and
$B=\{f(s_2,y_2),f(s_4,y_4),f(s_8,y_8),\newline f(s_6,y_6)\}$. In
particular, any couple $(s_i,s_{i+1})$ (resp. $(y_{(i)},y_{(i+1)})$)
contains
one $s$ (resp. $y$) of the group $A$ and one of $B$.

We can now come back to the analysis of the right-hand side of
(\ref{19}) and we can use the above detailed procedure to estimate
it by $2^{m-1} (s_0 x_0)^{\frac m2} m!\, (J_1+J_2)$, with
\begin{align*}
J_1 & =  n^m \int_{[0,T]^m\times [0,1]^m} \prod_{i_j\in\mathcal I}
\left( {\bf 1}_{[s_0,s_0']}(s_{i_j}){\bf 1}_{[x_0,x_0']}(y_{i_j})
f(s_{i_j},y_{i_j})^2 \right) \nonumber \\
&\quad \quad  \times \exp\left\{-n
x_0[(s_m-s_{m-1})+\cdots+(s_2-s_1)]\right\}
\nonumber \\
&\quad \quad  \times \exp\left\{-n
s_0[(y_{(m)}-y_{(m-1)})+\cdots+(y_{(2)}-y_{(1)})]\right\}
\nonumber \\
& \quad  \quad \times {\bf 1}_{\{ s_1\leq \dots\leq s_m\}} ds_1
\cdots ds_m dy_1\cdots dy_m.
\end{align*}
\begin{align*}
J_2 & =  n^m \int_{[0,T]^m\times [0,1]^m} \prod_{r_k\in\mathcal R}
\left( {\bf 1}_{[s_0,s_0']}(s_{r_k}){\bf 1}_{[x_0,x_0']}(y_{r_k})
f(s_{r_k},y_{r_k})^2 \right) \nonumber \\
&\quad \quad  \times \exp\left\{-n
x_0[(s_m-s_{m-1})+\cdots+(s_2-s_1)]\right\}
\nonumber \\
&\quad \quad  \times \exp\left\{-n
s_0[(y_{(m)}-y_{(m-1)})+\cdots+(y_{(2)}-y_{(1)})]\right\}
\nonumber \\
& \quad  \quad \times {\bf 1}_{\{ s_1\leq \dots\leq s_m\}} ds_1
\cdots ds_m dy_1\cdots dy_m.
\end{align*}
We will only deal with the term $J_1$, since $J_2$ can be treated
using exactly the same arguments. The idea is to integrate in $J_1$
with respect to $s_{r_k}, y_{r_k}$, with $r_k\in\mathcal R$,  for
$k=1,\dots,\frac m2$. Recall that the variables $s_{r_k}$ (resp.
$y_{r_k}$) have been chosen in such a way that they only appear once
in each couple $(s_i,s_{i+1})$ (resp. $(y_{(i)},y_{(i+1)})$).
Observe that we have, for any $k=1,\dots,\frac m2$,
$$\int_{s_0}^{s_0'} \exp\left\{ -nx_0(s_{r_k}-s_i)\right\}
{\bf 1}_{\{ s_i \leq s_{r_k}\}} ds_{r_k}\leq C \frac 1n$$ or
$$\int_{s_0}^{s_0'} \exp\left\{ -nx_0(s_{i+1}-s_{r_k} )\right\}
{\bf 1}_{\{ s_{r_k}\leq s_{i+1} \}} ds_{r_k}\leq C \frac 1n,$$ for
some $s_i$ and $s_{i+1}$, depending on which position occupies
$s_{r_k}$ in the corresponding couple. For the integrals with
respect to $y_{r_k}$ one obtains the same type of bound. Therefore,
\begin{align}
J_1 & \leq C_m \int_{[0,T]^{\frac m2}\times [0,1]^{\frac m2}}
\prod_{j=1}^{\frac m2} \left( {\bf 1}_{[s_0,s_0']}(s_{i_j}){\bf
1}_{[x_0,x_0']}(y_{i_j}) f(s_{i_j},y_{i_j})^2 \right) ds_{i_1}
\cdots ds_{i_{\frac m2}} dy_{i_1}\cdots
dy_{i_{\frac m2}} \nonumber\\
& = C_m \left( \int_0^T\int_0^1 {\bf 1}_{[s_0,s_0']}(s){\bf
1}_{[x_0,x_0']}(y) f(s,y)^2 \, dy ds\right)^{\frac m2}. \label{20}
\end{align}
As it has been mentioned, one can use the same arguments to get the
same upper bound for $J_2$. Hence, the right-hand side of (\ref{19})
can be estimated by (\ref{20}), and this concludes the proof.

\hfill $\Box$
\end{dem}

\medskip

\section{Convergence in law for the Donsker kernels}
\label{donsker}

In this section, we aim to prove that the hypothesis of Theorem
\ref{quasilinear} are satisfied in the case where the approximation
sequence is constructed in terms of the Donsker kernels. Namely, we
consider $\{Z_k,\,k\in\N^2\}$ an independent family of identically
distributed and centered random variables, with $E(Z_k^2)=1$ for all
$k\in\N^2$, and such that $E(|Z_k|^m)<+\infty$ for all $k\in\N^2$,
and some even number $m\ge 10$. Then, for all $n\geq 1$ and
$(t,x)\in[0,T]\times[0,1]$, we define the kernels
$$\theta_n(t,x)=n\sum_{k=(k^1,k^2)\in\N^2}Z_k\, {\bf 1}_{[k^1-1,k^1)\times[k^2-1,k^2)}(t n,x n).$$
Let us remind that the approximation sequence is given by 
\beq
X_n(t,x)=\int_0^t\int_0^1G_{t-s}(x,y)\theta_n(s,y)dyds,\; (t,x)\in
[0,T]\times [0,1].
\label{27}
\eeq 

Recall that Hypothesis \ref{hyp1} is a consequence of the extension of Donsker's theorem to the plane (see, for
instance, \cite{wichura}). On the other hand, we have the following result:


\begin{lema}\label{lema1} Let $\theta_n$ be the above defined Donsker kernels. Then,
there exists a positive constant $C_m$ such that, for any $f\in
L^2([0,T]\times[0,1])$, we have \begin{equation}E\left(\int_0^T
\int_0^1f(t,x)\theta_n(t,x)dx dt\right)^m\leq C_m\left(\int_0^T
\int_0^1 f^2(t,x) \,dx dt\right)^{\frac m2}, \label{donskerineq}
\end{equation}
 for all $n\geq 1$.
\end{lema}
\begin{obs}
Notice that, taking into account that $m\ge 10$, inequality
(\ref{donskerineq}) implies both Hypothesis \ref{hyp2} and
\ref{hyp3}, so that the hypotheses of Theorem \ref{quasilinear} are satisfied for the Donsker kernels.
\end{obs}

{\it{Proof of Lemma \ref{lema1}:}}

First, we observe that we can write
\begin{eqnarray}
&&E\left(\int_0^T \int_0^1f(t,x)\theta_n(t,x)dx dt\right)^m\label{equ1}\\
&=&\int_{[0,T]^m\times[0,1]^m}f(t_1,x_1)\cdots
f(t_m,x_m)E\left[\prod_{j=1}^m\theta_n(t_j,x_j)\right]dt_1\cdots
dt_m dx_1\cdots dx_m\nonumber
\end{eqnarray}
By definition of $\theta_n$,
\begin{eqnarray*}
&&E\left[\prod_{j=1}^m\theta_n(t_j,x_j)\right]\\&=&n^{m}E\left[\prod_{j=1}^m\left(\sum_{k=(k^1,k^2)\in\N^2}Z_k\,
{\bf 1}_{[k^1-1,k^1)}(t_j n ){\bf 1}_{[k^2-1,k^2)}(x_j
n )\right)\right]\\
&=&n^{m}\sum_{k_1,\dots,k_m\in\N^2}E(Z_{k_1}\cdots
Z_{k_m})\prod_{j=1}^m\left({\bf 1}_{[k_j^1-1,k_j^1)}(t_j n ){\bf
1}_{[k_j^2-1,k_j^2)}(x_j n )\right).
\end{eqnarray*}
Notice that, by hypothesis,  $E(Z_{k_1}\cdots Z_{k_m})=0$ if, for
some $j\in\{1,\dots,m\}$, we have that $k_j\neq k_l$ for all
$l\in\{1,\dots,m\}\setminus\{j\}$; that is, if some variable
$Z_{k_j}$ appears only once in the product $Z_{k_1}\cdots Z_{k_m}$.

On the other hand, since $E(|Z_k|^m)<\infty$ for all $k\in\N^2$, then $E(Z_{k_1}\cdots Z_{k_m})$ is bounded for all
$k_1,\dots,k_m\in\N^2$ . Hence,
$$E\left[\prod_{j=1}^m\theta_n(t_j,x_j)\right]\leq n^{m}C_m\sum_{(k_1,\dots,k_m)\in
A^m}\prod_{j=1}^m\left({\bf 1}_{[k_j^1-1,k_j^1)}(t_j n ){\bf
1}_{[k_j^2-1,k_j^2)}(x_j n )\right),$$ with
$$A^m=\{(k_1,\dots,k_m)\in\N^{2m};\, \textrm{ for all }l\in\{1,\dots,m\}, k_l=k_j \textrm{ for some }
j\in\{1,\dots,m\}\setminus\{l\}\}.$$ Notice that we have the
following estimation:
$$\sum_{(k_1,\dots,k_m)\in
A^m}\prod_{j=1}^m\left({\bf 1}_{[k_j^1-1,k_j^1)}(t_j n ){\bf
1}_{[k_j^2-1,k_j^2)}(x_j n )\right)\leq {\bf
1}_{D^m}(t_1,\dots,t_m;x_1,\dots,x_m),$$ where $D^m$ denotes the set
of $(t_1,\dots,t_m;x_1,\dots,x_m)\in[0,T]^{m}\times[0,1]^m$
satisfying the following property: for all $l\in\{1,\dots,m\}$,
there exists $j\in\{1,\dots,m\}\setminus\{l\}$ such that
$|t_j-t_l|<\frac1n$ and $|x_j-x_l|<\frac1n$ and, moreover, if there
is some $r\neq j,l$ verifying $|t_l-t_r|< \frac 1n$ and
$|x_l-x_r|<\frac 1n$, then $|t_j-t_r|<\frac 1n$ and $|x_j-x_r|<\frac
1n$.

Next, observe that we can bound
$I_{D^m}(t_1,\dots,t_m;x_1,\dots,x_m)$ by a finite sum of products
of indicators, where in each product of indicators there appear all
the $m$ variables $t_1,\dots,t_m$ and all the $m$ variables
$x_1,\dots,x_m$, but each indicator concerns only two or three of
them. Moreover, each variable only appears in one of the indicators
of each product and, whenever we have some indicator concerning two
variables $t_j$ and $t_l$, (respectively three variables $t_j$,
$t_l$ and $t_r$), we have the same indicator for the variables $x_j$
and $x_l$, (respectively for the variables $x_j$, $x_l$ and $x_r$).
Therefore, expression (\ref{equ1}) can be bounded by a finite sum of
products of the following two kinds of terms:
\begin{itemize}
\item[(i)] For some  $l,j\in\{1,\dots,m\}$ such that $l\neq j$,
\beq C_mn^2\int_{[0,T]^{2}\times[0,1]^2}|f(t_l,x_l)||f(t_j,x_j)|{\bf
1}_{[0,\frac1n)}(|t_j-t_l|){\bf
1}_{[0,\frac1n)}(|x_j-x_l|)dt_jdt_ldx_jdx_l. \label{26} \eeq
\item[(ii)] For some $l,j,r\in\{1,\dots,m\}$ such that $l\neq j$, $l\neq
r$ and $r\neq j$,
\begin{eqnarray*}
&&C_mn^{3}\int_{[0,T]^{3}\times[0,1]^3}|f(t_l,x_l)||f(t_j,x_j)||f(t_r,x_r)|\\&&
\quad \times {\bf 1}_{[0,\frac1n)}(|t_j-t_l|){\bf
1}_{[0,\frac1n)}(|t_l-t_r|){\bf 1}_{[0,\frac1n)}(|t_j-t_r|)
\\&&\quad \times
{\bf 1}_{[0,\frac1n)}(|x_j-x_l|){\bf 1}_{[0,\frac1n)}(|x_l-x_r|){\bf
1}_{[0,\frac1n)}(|x_j-x_r|)dt_jdt_ldt_rdx_jdx_ldx_r.
\end{eqnarray*}
\end{itemize}
Then, it turns out that, in order to conclude the proof, it suffices
to bound the first type of term (i) by $C_m\int_0^T
\int_0^1f^2(t,x)\,dx dt$ and the second one (ii) by
$C_m\left(\int_0^T \int_0^1f^2(t,x)\, dx dt\right)^{\frac32}$.

Let us use first the fact that, for all $a,b\in\R$, $2ab\leq
a^2+b^2$, so that a term of the form (\ref{26}) can be bounded, up
to some constant, by
\begin{eqnarray*}
&&C_mn^2\int_{[0,T]^{2}\times[0,1]^2}f^2(t_l,x_l) {\bf 1}_{[0,\frac1n)}(|t_j-t_l|){\bf 1}_{[0,\frac1n)}(|x_j-x_l|)dt_jdt_ldx_jdx_l\\
&\leq&C_m\int_0^T \int_0^1f^2(t,x) \, dx dt.
\end{eqnarray*}
On the other hand, using that for all $a,b,c\in\R^{+}$,
$2abc\leq(ab^2+ac^2)$, we can study the terms of type $(ii)$ in the
following way:
\begin{eqnarray*}
&&C_mn^{3}\int_{[0,T]^{3}\times[0,1]^3}|f(t_l,x_l)||f(t_j,x_j)||f(t_r,x_r)|\\&&\times
{\bf 1}_{[0,\frac1n)}(|t_j-t_l|){\bf 1}_{[0,\frac1n)}(|t_l-t_r|){\bf
1}_{[0,\frac1n)}(|t_j-t_r|)
\\&&\times
{\bf 1}_{[0,\frac1n)}(|x_j-x_l|){\bf 1}_{[0,\frac1n)}(|x_l-x_r|){\bf
1}_{[0,\frac1n)}(|x_j-x_r|)dt_jdt_ldt_rdx_jdx_ldx_r
\\&\leq&C_mn^{3}\int_{[0,T]^{3}\times[0,1]^3}|f(t_l,x_l)|f^2(t_j,x_j) \\&&\times
{\bf 1}_{[0,\frac1n)}(|t_j-t_l|){\bf 1}_{[0,\frac1n)}(|t_l-t_r|){\bf
1}_{[0,\frac1n)}(|t_j-t_r|)
\\&&\times
{\bf 1}_{[0,\frac1n)}(|x_j-x_l|){\bf 1}_{[0,\frac1n)}(|x_l-x_r|){\bf
1}_{[0,\frac1n)}(|x_j-x_r|)dt_jdt_ldt_rdx_jdx_ldx_r
\\
&\leq&C_mn\int_{[0,T]^{2}\times[0,1]^2}|f(t_l,x_l)|f^2(t_j,x_j)\,
{\bf 1}_{[0,\frac1n)}(|t_j-t_l|) {\bf
1}_{[0,\frac1n)}(|x_j-x_l|)dt_jdt_ldx_jdx_l
\\&=&C_mn\int_0^T \int_0^1|f(t_l,x_l)|\left(\int_0^T
\int_0^1f^2(t_j,x_j)\, {\bf 1}_{[0,\frac1n)}(|t_j-t_l|) {\bf
1}_{[0,\frac1n)}(|x_j-x_l|)dt_jdx_j\right)dt_ldx_l.
\end{eqnarray*}
At this point, we apply Cauchy-Schwarz inequality, so that the
latter expression can be estimated by
\begin{eqnarray*}
&&C_mn\left(\int_0^T \int_0^1f^2(t_l,x_l)\,
dt_ldx_l\right)^{\frac12}\\&& \times \left(\int_0^T
\int_0^1\left(\int_0^T \int_0^1f^2(t_j,x_j) {\bf
1}_{[0,\frac1n)}(|t_j-t_l|)
{\bf 1}_{[0,\frac1n)}(|x_j-x_l|)dt_jdx_j\right)^2dt_ldx_l\right)^{\frac12}\\
&=&C_mn\left(\int_0^T \int_0^1f^2(t_l,x_l)\,
dt_ldx_l\right)^{\frac12}\\&& \times
\left(\int_{[0,T]^3\times[0,1]^3}f^2(t_j,x_j) f(t_p,x_p)^2 {\bf
1}_{[0,\frac1n)}(|t_j-t_l|) {\bf
1}_{[0,\frac1n)}(|x_j-x_l|)\right.\\&&\left.\times {\bf
1}_{[0,\frac1n)}(|t_p-t_l|) {\bf 1}_{[0,\frac1n)}(|x_p-x_l|)
dt_jdt_pdt_ldx_jdx_pdx_l\right)^{\frac12}\\
&\leq&C_m\left(\int_0^T \int_0^1f^2(t,x)dx dt\right)^{\frac32}.
\end{eqnarray*}

This finishes the proof of the lemma. \hfill $\Box$

\vspace{0.3cm}

\appendix
\section{Appendix}
In this appendix, we give a sketch of the proof of Lemma \ref{lema5}
and discuss the relation between our results and those of
Manthey in \cite{manthey1} (see also \cite{manthey2}).

\bigskip

{\it Proof of Lemma \ref{lema5}:} As we have already
pointed out, we will only give the main lines of the proof.

Let $f\in E$ and $h\in \C^1(\R)$ having a bounded derivative. We aim
to prove that, for any $\eta>0$, it holds \beq \left|
E[h(J^n(f))]-E[h(J(f))]\right| < \eta,\label{29}\eeq for
sufficiently big $n$. For this, the idea is to consider an element
$g$ in $D$ which is close to $f$ with respect to the norm
$\|\cdot\|$. Then, one splits the left-hand side of (\ref{29}) in
several terms, which can be easily treated using the following facts:
\begin{enumerate}
  \item When $f$ is replaced by $g$, we have that the left-hand side of (\ref{29}) converges to zero, by hypothesis.
  \item One keeps control of the remaining terms using that $h$
  defines a Lipschitz function and that (\ref{30}) and (\ref{30'}) hold.
\end{enumerate}

\hfill $\Box$

\bigskip

\noindent{\bf Relation with Manthey results}

\medskip

\noindent In \cite{manthey1}, the author considers the family of
processes  $\{X_n,\, n\in\mathbb N\}$ such that each $X_n$ is the
mild solution of the equation
\begin{equation*}
\frac{\partial X_n}{\partial t}(t,x)-\frac{\partial ^2 X_n}{\partial
x^2} (t,x)=\theta_n(t,x),\quad (t,x)\in [0,T]\times [0,1],
\end{equation*}
with null initial condition and Dirichlet boundary conditions. The processes $\theta_n$ are correlated noises satisfying the
following conditions:
\begin{enumerate}[(i)]
\item For all $(t,x)\in [0,T]\times [0,1]$, $$\int_0^t\int_0^x \theta_n(s,y)^2 dyds<\infty,\;a.s.$$
\item For each $ m\in \N$ and $(t_1,x_1),\dots,(t_m,x_m)\in [0,T]\times [0,1]$,  the random vector
$$\left(\int_0^{t_1}\int_0^{x_1} \theta_n(s,y)dyds,\ldots ,\int_0^{t_m}\int_0^{x_m}
\theta_n(s,y)dy\right)$$ converges weakly to
$(W(t_1,x_1),\ldots,W(t_m,x_m))$, where we recall that $\{W(t,x),\; (t,x)\newline \in [0,T]\times [0,1]\}$ denotes a Brownian sheet.
\item For all $(t,x)\in [0,T]\times [0,1]$, $E[\theta_n(t,x)]=0$.
\item There exists $n_0\in\mathbb N$ such that
$$\sup_{ \bistackrel{n\ge n_0}{(t,x)\in [0,T]\times [0,1]}{} }\int_0^{t}\int_0^{x}\left|E\big[\theta_n(s,y)\theta_n(t,x)\big]\right|dyds<\infty.$$
\end{enumerate}
Under these hypotheses, it has been proved that $X_n$ converges weakly, in
the sense of the convergence of finite dimensional distributions, to
the process $X$ which is the mild solution of
\begin{equation*}
\frac{\partial X}{\partial t}(t,x)-\frac{\partial ^2 X}{\partial
x^2} (t,x)=\dot W(t,x),\quad (t,x)\in [0,T]\times [0,1].
\end{equation*}
Furthermore, it is showed that, if the processes $\theta_n$ are Gaussian, the
convergence also holds in $\mathcal C([0,T]\times [0,1])$. These
results are extended to the quasi-linear equation (\ref{mild})
in \cite{manthey2}.

First, it is worth pointing out that one can easily see that condition (iii) is not essential in the proof. Moreover,
in Manthey's result, condition (iv) is stated in a weaker form, though
we believe that, in his proof, the statement which has been used is indeed condition (iv) as stated above
(see the last inequality in p. 163 of \cite{manthey1}).

Secondly, one can easily see that condition (iv) stated above implies
Hypothesis \ref{hyp2}  with $q=2$. Therefore, the hypotheses
assumed in Proposition \ref{dimfin} (which assures the convergence of
the finite dimensional distributions) are weaker than (i)-(iv).

Eventually, processes $\theta_n$ given by the Kac-Stroock
processes and the Donsker kernels are not Gaussian so that, if
conditions (i)-(iv) were satisfied, using Manthey's result only
convergence of the finite dimensional distributions could be
obtained. In fact, it is straightforward to check that the Donsker kernels
satisfy these conditions, but condition (iv) fails for the
Kac-Stroock processes. This is proved in the following lemma:

\begin{lema}\label{fails}
Assume that $\{\theta_n(t,x),\; (t,x)\in
[0,T]\times[0,1]\}$, $n\geq 1$, is the Kac-Stroock process
(\ref{7}). Then, the family $\{\theta_n,\,n\in \mathbb N\}$ does not
satisfy condition (iv) above.
\end{lema}
\begin{dem} We will show that, when
$$\theta_n(s, y)=n\sqrt{sy}\,(-1)^{N_n(s,y)},$$
then the quantity
$$\int_0^T\int_0^1 \left|E[\theta_n(s,y)\theta_n(t,x)]\right|\,dy ds$$ is not
uniformly bounded in $n,\,t$ and $x$. Indeed, it holds
\beq\label{mant}
\int_0^T\int_0^1 |E(\theta_n(s,y)\theta_n(t,x))|\,dy ds=\intdob
n^2\sqrt{sytx} \, E\left[(-1)^{N_n(s,y)+N_n(t,x)}\right]dy ds.
\eeq
Owing to the proof of \cite[Lemma 3.1]{bjt} (see also \cite[Lemma 3.2]{bj2000}), we have that
\begin{eqnarray*}
E\left[(-1)^{N_n(s,y)+N_n(t,x)}\right]\,\,= & e^{-2n[(t-s)x+(x-y)s]}\,\mathbf{1}_{\{s\le
t,y\le x\}}\,+\, e^{-2n[(t-s)x+(y-x)s]}\,\mathbf{1}_{\{ s\le t,y\ge
x\}}\\
 & + e^{-2n[(s-t)y+(x-y)t]}\,\mathbf{1}_{\{ s\ge t, y\le x\}}\,+\,
e^{-2n[(s-t)y+(y-x)t]}\,\mathbf{1}_{\{s\ge t, y\ge x\}}.
\end{eqnarray*}
 Then,
expression (\ref{mant}) is the sum of four positive integrals. It is clear that one
of them is given by
$$I(n,t,x)=\int_0^t\int_x^1 n^2 \sqrt{sytx} \, e^{-2n[(t-s)x+(y-x)s]} dy
ds.$$ We will check that this integral is not uniformly bounded. In
fact,
\begin{eqnarray}
\sup_{n,\,t,\,x}I(n,t,x)&\ge&\sup_n
I\left(n,T,\frac1n\right)\nonumber \\ & = &\sup_n \sqrt{T} \int_0^T\int_{\frac1n}^1
n^2\frac{\sqrt{sy}}{\sqrt n}\,e^{-2n[(T-s)\frac1n+(y-\frac1n)s]} dy
ds\nonumber\\&=&\sqrt{T} e^{-2T}\sup_n \int_0^T\int_{\frac1n}^1 n^{\frac32}
\sqrt{sy}\,e^{4s}e^{-2nys}  dy ds\nonumber\\&=&\sqrt{T} e^{-2T}\sup_n\int_0^T\int_1^n
\sqrt{sz}\,e^{4s}e^{-2zs}  dz ds\nonumber \\ & = & \sqrt{T} e^{-2T}\int_0^T\int_1^{\infty}
\sqrt{sz}\,e^{4s}e^{-2zs}  dz ds.
\label{int}
\end{eqnarray}
Let us apply the change of coordinates $v=sz$, for any fixed $s$, and then Fubini Theorem in the last integral of (\ref{int}), so that we end up with
$$\sup_{n,\,t,\,x}I(n,t,x) \geq \sqrt{T} e^{-2T} \int_0^{+\infty} \sqrt{v}\, e^{-2v} \left(\int_0^{v\land T} \frac 1s \, e^{4s} ds\right) dv,$$
and the latter is clearly divergent. This fact concludes the proof. \hfill $\Box$
\end{dem}

\bigskip





\end{document}